\documentclass[11pt]{article}
\usepackage[T1]{fontenc}
\usepackage[utf8]{inputenc}
\usepackage{float}
\usepackage[a4paper,margin=1in]{geometry}
\usepackage{amsmath,amssymb,amsthm,mathtools}
\usepackage{microtype}
\usepackage[hidelinks]{hyperref}
\usepackage{graphicx,tikz}
\usetikzlibrary{arrows.meta,calc}

\usepackage[T1]{fontenc}
\usepackage{lmodern,amsmath,amssymb}
\usetikzlibrary{arrows.meta,calc}
\newcommand{\R}{\mathbb R}
\newcommand{\C}{\mathcal C}
\newcommand{\M}{\mathcal M}
\newtheorem{theorem}{Theorem}[section]
\newtheorem{proposition}[theorem]{Proposition}
\newtheorem{lemma}[theorem]{Lemma}

\theoremstyle{definition}
\newtheorem{definition}[theorem]{Definition}
\newtheorem{example}[theorem]{Example}
\theoremstyle{remark}
\newtheorem{remark}[theorem]{Remark}
\numberwithin{equation}{section}

\title{Supercritical Nodal Solutions on Annuli:
Symmetry Breaking and Multiplicity}
\author{ Alireza Khatib \thanks{Corresponding author. Universidade Federal do Amazonas, Manaus-AM, Brazil, alireza@ufam.edu.br
 }  \quad Olimpio H. Miyagaki\thanks{ Universidade Federal de São Carlos, São Carlos, Brazil, olimpio@ufscar.br } \quad Somayeh Mousavinasr \thanks{Universidade Federal do Amazonas, Manaus-AM, Brazil, somayeh@ufam.edu.br
 } \quad 
}
\date{}

\begin{document}
\maketitle
\begin{abstract}
We establish existence and multiplicity results for nonradial
sign-changing solutions of the Dirichlet problem
$-\Delta u=a(|x|)|u|^{p-2}u$ on annuli in $\mathbb R^N$,
where $N\geq4$, $a$ is a positive radial coefficient,
and $p$ may exceed the Sobolev critical exponent.
For each decomposition $N=m+n$, with $m\geq n\geq2$,
we minimize the energy on a nodal Nehari set within
a signed cone of angularly nondecreasing,
$O(m)\times O(n)$-invariant functions.
Whenever $2<p<2(m+1)/(m-1)$, the minimum is attained,
and every minimizer is a classical solution with exactly
two nodal domains.
A cone-preserving auxiliary Dirichlet operator and
a two-parameter energy deformation recover the full
equation from the constrained minimization.
We prove that these minimizers are nonradial if
$(p-2)\beta_\Omega>2N$, where $\beta_\Omega$ is
the optimal Dirichlet Hardy constant associated with
the weight $|x|^{-2}$.
Under a common admissible exponent range satisfying
this condition, different orthogonal decompositions
yield distinct solution pairs $\{\pm u\}$.
An explicit lower bound for $\beta_\Omega$ provides
supercritical regimes on sufficiently thin annuli
in which the nonradiality and multiplicity results apply.
\end{abstract}

\section{Introduction}
\label{sec:introduction}

We study the Dirichlet problem
\begin{equation}\label{eq:problem}
\begin{cases}
-\Delta u=a(|x|)|u|^{p-2}u &\text{in }\Omega,\\
u=0 &\text{on }\partial\Omega,
\end{cases}
\end{equation}
on the annulus
\[
\Omega=A_{R_-,R_+}
:=\{x\in\mathbb R^N:R_-<|x|<R_+\},
\qquad 0<R_-<R_+,\quad N\geq4.
\]
Throughout the paper, the coefficient satisfies
\begin{equation}\label{eq:a}
a\in C^1([R_-,R_+]),
\qquad 0<a_0\leq a(r)\leq a_1.
\end{equation}
Our aim is to construct nonradial sign-changing solutions with
exactly two nodal domains, for powers that may exceed the Sobolev
critical exponent
\[
2_N^*:=\frac{2N}{N-2}.
\]
We obtain these solutions by minimizing the energy under a signed
angular-order constraint and derive an explicit condition relating
symmetry breaking to the geometry of the annulus.

The variational study of nodal solutions connects energy minimization
with the geometry of the positive and negative regions; see also
\cite{BonheureEtAl2022}.
In the subcritical setting, Bartsch, Weth, and Willem
\cite{BartschWethWillem2005} established partial symmetry properties
of least-energy nodal solutions by polarization. Further symmetry
results for least-energy positive and nodal solutions were obtained
by Van Schaftingen and Willem \cite{VanSchaftingenWillem2008}.
Bartsch and
Degiovanni \cite{BartschDegiovanni2006} proved nonradiality of
least-energy nodal solutions on balls and annuli under subcritical
growth assumptions. These results show that radial symmetry of the
domain need not be inherited by nodal minimizers. Nehari-type
minimization and its relation to nodal partitions were developed
further by Conti, Terracini, and Verzini
\cite{ContiTerraciniVerzini2002}. Variational bounds on the number
of nodal domains were established by Bartsch and Weth
\cite{BartschWeth2003}.

For positive solutions on annuli, Gladiali, Grossi, Pacella, and
Srikanth \cite{GladialiGrossiPacellaSrikanth2011} proved nonradial
bifurcation from radial solutions. For sign-changing solutions,
Amadori and Gladiali \cite{AmadoriGladiali} constructed
global continua of nonradial sign-changing solutions to the
Lane--Emden problem through bifurcation and fixed-point index
arguments. In the supercritical setting, Clapp and Pacella
\cite{ClappPacella2017} established the existence of infinitely many
nonradial nodal solutions on annuli in suitable dimensions and
exponent ranges, using prescribed symmetries. They also studied
asymptotic profiles as the exponent approaches the critical value
and the hole shrinks. The present work addresses solutions with
\emph{angular order and exactly two nodal domains}, together with
a quantitative criterion for nonradiality and a multiplicity result
associated with different orthogonal decompositions.

Symmetry provides a way to treat powers above $2_N^*$ by improving
Sobolev embeddings; see Hebey and Vaugon
\cite{HebeyVaugon1997} for the general role of group orbits in
such improvements. Cones defined by symmetry and angular monotonicity
also play an important role in supercritical elliptic problems;
see, for instance,
\cite{BoscagginEtAl2023,BoscagginEtAl2024,BoscagginEtAl2026}.
Cowan and Moameni \cite{CowanMoameni,CowanMoameni2024} developed
cone-based variational methods for supercritical elliptic problems,
including existence, multiplicity, and symmetry breaking of positive
solutions. A general critical point framework on convex subsets
was established by Moameni \cite{Moameni2020}. Related methods for
positive solutions of Neumann problems on annuli, with supercritical
power and exponential nonlinearities, were developed by Khatib and
Mousavinasr \cite{KhatibMousavinasr2026Neumann}.
Here the exponent range is chosen so that compactness already holds
in the full invariant Sobolev space. The angular constraint is used
to preserve order under the auxiliary Dirichlet operator and to
control the nodal-domain count. In particular, we impose angular
monotonicity on the signed function itself: its positive and negative
parts then have opposite angular orders.

To state our results, fix a decomposition
\[
N=m+n,\qquad m\geq n\geq2,
\qquad G=O(m)\times O(n).
\]
Write $x=(x',x'')\in\mathbb R^m\times\mathbb R^n$ and introduce
the orbit coordinates
\[
|x'|=r\cos\theta,\qquad |x''|=r\sin\theta,
\qquad 0\leq\theta\leq\frac\pi2.
\]
Let $H_{0,G}^1(\Omega)$ denote the space of $G$-invariant functions
in $H_0^1(\Omega)$, and set
\[
Q_{m,n}:=\frac{2(m+1)}{m-1},
\qquad
\mathcal C_{m,n}
:=\{u\in H_{0,G}^1(\Omega):u_\theta\geq0
\text{ a.e.}\}.
\]
The embedding of $H_{0,G}^1(\Omega)$ into $L^p(\Omega)$ is compact
for $2<p<Q_{m,n}$. Since $n\geq2$, we have $Q_{m,n}>2_N^*$,
so this range includes supercritical powers.

For $2<p<Q_{m,n}$, the energy
\[
J(u):=\frac12\int_\Omega|\nabla u|^2\,dx
-\frac1p\int_\Omega a(|x|)|u|^p\,dx
\]
is of class $C^1$ on $H_{0,G}^1(\Omega)$. With
$u^+:=\max\{u,0\}$ and $u^-:=\max\{-u,0\}$, define the nodal
Nehari set
\[
\begin{aligned}
\mathcal M_{m,n}:=\biggl\{u\in\mathcal C_{m,n}:\;&
u^+\not\equiv0,\quad u^-\not\equiv0,\\
&\int_\Omega|\nabla u^\pm|^2\,dx
=\int_\Omega a(|x|)(u^\pm)^p\,dx\biggr\}.
\end{aligned}
\]
Each sign therefore satisfies its own Nehari identity. We also
introduce the annular Hardy constant
\begin{equation}\label{eq:Dirichlet-Hardy-constant}
\beta_\Omega:=
\inf_{v\in H_0^1(\Omega)\setminus\{0\}}
\frac{\displaystyle\int_\Omega|\nabla v|^2\,dx}
{\displaystyle\int_\Omega |x|^{-2}v^2\,dx}>0.
\end{equation}

\begin{theorem}[Existence and nonradiality]\label{thm:main}
Assume \eqref{eq:a}, $N=m+n$, $m\geq n\geq2$, and
$2<p<Q_{m,n}$. The infimum of $J$ on $\mathcal M_{m,n}$
is positive and attained. Every minimizer is a weak solution of
\eqref{eq:problem}, belongs to $C^{2,\alpha}(\overline\Omega)$
for some $\alpha\in(0,1)$, and has exactly two nodal domains.
If
\begin{equation}\label{eq:nodal-Hardy-condition}
(p-2)\beta_\Omega>2N,
\end{equation}
then every such minimizer is nonradial.
\end{theorem}

The energy characterization in this theorem is relative to the
prescribed angular cone. The two nodal domains are obtained as part
of the construction, without fixing their common interface.
Condition~\eqref{eq:nodal-Hardy-condition} gives a sufficient
criterion for symmetry breaking in terms of the exponent and the
annular Hardy constant.

The same constant can be used for several orthogonal decompositions.
In the next theorem, $O(N-j)\times O(j)$ acts on the first $N-j$
and last $j$ standard coordinates.

\begin{theorem}[Multiplicity]\label{thm:main-multiplicity}
Assume \eqref{eq:a}. Let $\mathcal J$ be a nonempty subset of
$\{2,\ldots,\lfloor N/2\rfloor\}$ and suppose
\begin{equation}\label{eq:multiplicity-range}
2+\frac{2N}{\beta_\Omega}<p<
\min_{j\in\mathcal J}\frac{2(N-j+1)}{N-j-1}.
\end{equation}
Then \eqref{eq:problem} has at least $|\mathcal J|$ distinct
pairs $\{\pm u_j\}$ of nonradial sign-changing solutions.
For each $j\in\mathcal J$, the solution $u_j$ is
$O(N-j)\times O(j)$-invariant and has exactly two nodal domains.
Moreover, $u_j\neq\pm u_\ell$ for $j\neq\ell$.
\end{theorem}
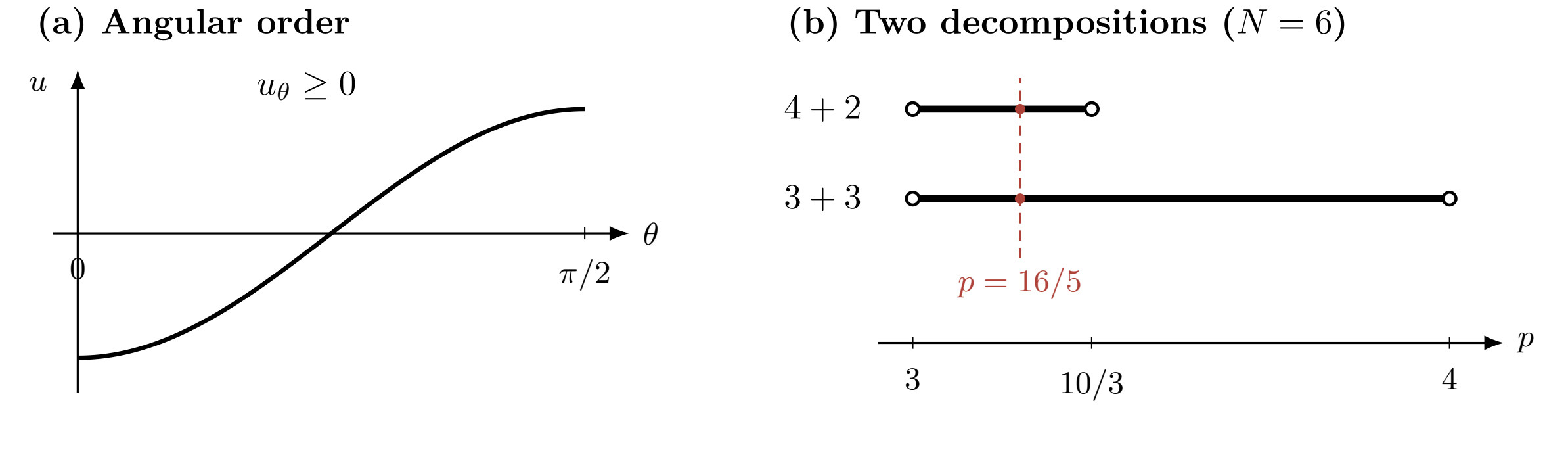
\begin{figure}[H]
  \centering
  \resizebox{\textwidth}{!}{
\definecolor{nodalred}{HTML}{AF4238}
\begin{tikzpicture}[x=1cm,y=1cm,>=Latex,font=\small,text=black,
  axis/.style={->,line width=.6pt,draw=black},
  paneltitle/.style={anchor=west,font=\small\bfseries}]
  \path[use as bounding box] (0,0) rectangle (15.65,4.55);
  \node[paneltitle] at (.15,4.3) {(a) Angular order};
  \draw[axis] (.7,.6)--(.7,3.85);
  \draw[axis] (.45,2.2)--(6.25,2.2);
  \node[anchor=east,font=\footnotesize] at (.53,3.7) {$u$};
  \node[right,font=\footnotesize] at (6.25,2.2) {$\theta$};
  \draw[black,line width=1.25pt]
    plot[domain=0:180,samples=101]
      ({.7+5.1*\x/180},{2.2-1.25*cos(\x)});
  \node[font=\small] at (3.0,3.68) {$u_\theta\geq0$};
  \draw[line width=.4pt] (.7,2.14)--(.7,2.26);
  \draw[line width=.4pt] (5.8,2.14)--(5.8,2.26);
  \node[below,font=\footnotesize] at (.7,2.08) {$0$};
  \node[below,font=\footnotesize] at (5.8,2.08) {$\pi/2$};

  \node[paneltitle] at (7.7,4.3) {(b) Two decompositions ($N=6$)};
  \node[anchor=east] at (8.72,3.45) {$4+2$};
  \node[anchor=east] at (8.72,2.55) {$3+3$};
  \draw[black,line width=2pt] (9.1,3.45)--(10.9,3.45);
  \draw[black,line width=.85pt,fill=white] (9.1,3.45) circle (.065);
  \draw[black,line width=.85pt,fill=white] (10.9,3.45) circle (.065);
  \draw[black,line width=2pt] (9.1,2.55)--(14.5,2.55);
  \draw[black,line width=.85pt,fill=white] (9.1,2.55) circle (.065);
  \draw[black,line width=.85pt,fill=white] (14.5,2.55) circle (.065);
  \draw[nodalred,densely dashed,line width=.6pt] (10.18,1.95)--(10.18,3.76);
  \fill[nodalred] (10.18,3.45) circle (1.5pt);
  \fill[nodalred] (10.18,2.55) circle (1.5pt);
  \node[font=\footnotesize,text=nodalred] at (10.18,1.7) {$p=16/5$};
  \draw[axis] (8.75,1.1)--(15.05,1.1);
  \node[right,font=\footnotesize] at (15.05,1.1) {$p$};
  \foreach \xx in {9.1,10.9,14.5}
    \draw[line width=.4pt] (\xx,1.04)--(\xx,1.16);
  \node[below,font=\footnotesize] at (9.1,.97) {$3$};
  \node[below,font=\footnotesize] at (10.9,.97) {$10/3$};
  \node[below,font=\footnotesize] at (14.5,.97) {$4$};
\end{tikzpicture}
}
\caption{(a) A schematic signed angular profile with
$u_\theta\geq0$.
(b) Supercritical admissible exponent ranges for
$6=4+2$ and $6=3+3$ on $A_{1,11/10}$, with
$2_6^*=3$. The dashed line marks $p=16/5$.
The Hardy condition holds throughout both displayed
intervals.}

  \label{fig:signed-cone-multiplicity}
\end{figure}
The conditions above can be checked directly on thin annuli.
Writing $d:=R_+-R_-$, we prove
\[
\beta_\Omega\geq
\frac{\pi^2R_-^2}{d^2}
\left(\frac{R_-}{R_+}\right)^{N-1}.
\]
For fixed $R_->0$, the right-hand side tends to infinity as
$d\downarrow0$. Thus, for every fixed exponent in the relevant
compactness range, the nonradiality condition holds when the annulus
is sufficiently thin. A concrete supercritical example is illustrated in
Figure~\ref{fig:signed-cone-multiplicity} and verified
in Example~\ref{ex:two-pairs}.

The main variational step is to show that a minimizer on
$\mathcal M_{m,n}$ satisfies the full equation. We associate with
each invariant function $u$ the auxiliary Dirichlet solution $R(u)$
with source $a(|x|)|u|^{p-2}u$. This operator preserves the signed
angular cone. Convexity of the nonlinear potential gives an explicit
energy decrease along the segment from $u$ to $R(u)$. Applying this
decrease to a two-parameter family and using the
Poincar\'e--Miranda theorem shows that a nodal minimizer must satisfy
$u=R(u)$. The argument adapts the deformation strategy used in nodal
Nehari minimization; see, for example,
\cite[Lemma~5.3]{ContiTerraciniVerzini2002}. Here the deformation
must also preserve the signed angular constraint.

Once the full equation is recovered, angular order ensures that
restriction to an individual positive or negative nodal component
preserves the corresponding cone condition. Minimality then excludes
additional nodal domains. Finally, angular perturbations of radial
minimizers yield \eqref{eq:nodal-Hardy-condition}, and comparison
of the prescribed symmetry classes gives the multiplicity result.

The paper is organized as follows.
Section~\ref{sec:setting} establishes the invariant Sobolev embeddings
and the properties of the signed angular cone.
Section~\ref{sec:minimum} proves attainment of the nodal level.
Section~\ref{sec:criticality} recovers the full weak equation through
the auxiliary operator and the cone-preserving deformation.
Section~\ref{sec:nodal-domains} establishes regularity and the exact
nodal-domain count. Section~\ref{sec:nonradiality} proves symmetry
breaking and multiplicity and gives the thin-annulus criterion.

\section{Functional setting and the signed angular cone}
\label{sec:setting}

Symmetry provides the compact embedding needed for minimization.
We then impose angular order through a signed cone and describe
the resulting orders of the positive and negative parts.
Finally, we explain how averaging over the symmetry group
allows the weak equation to be tested against every function
in $H_0^1(\Omega)$.

Fix $N=m+n$ with $m\geq n\geq2$, and let
\[
G:=O(m)\times O(n),\qquad
(A,B)(x',x''):=(Ax',Bx'').
\]
Set
\[
H:=H_{0,G}^1(\Omega),\qquad
(u,v)_H:=\int_\Omega\nabla u\cdot\nabla v\,dx,
\qquad \|u\|^2:=(u,u)_H.
\]
Here $H_{0,G}^1(\Omega)$ is the closed subspace of
$G$-invariant functions in $H_0^1(\Omega)$. We use the
Dirichlet norm on $H_0^1(\Omega)$ and write $\|u\|_q$ for
the $L^q(\Omega)$ norm.

An invariant function has an orbit profile $u(s,t)$, where
$s=|x'|$, $t=|x''|$, on
\[
\widehat\Omega:=\{(s,t)\in(0,\infty)^2:
R_-^2<s^2+t^2<R_+^2\}.
\]
In orbit polar coordinates,
\[
s=r\cos\theta,\qquad t=r\sin\theta,
\qquad u_\theta=s u_t-t u_s.
\]
Angular derivatives of Sobolev functions are understood weakly
on the interior of $\widehat\Omega$.
The orbit of $x=(x',x'')$ is
\[
Gx=|x'|\mathbb S^{m-1}\times |x''|\mathbb S^{n-1}.
\]
It is connected because $m,n\geq2$, with a zero-radius factor
interpreted as a singleton. This fact will be used in the
nodal-domain argument.

We first record an embedding for the whole invariant space $H$.
It supplies the compactness used in the following sections.

\begin{proposition}\label{prop:compactness}
Set $Q=Q_{m,n}=2(m+1)/(m-1)$. The embedding
$H\hookrightarrow L^q(\Omega)$ is continuous for $1\leq q\leq Q$
and compact for $1\leq q<Q$.
\end{proposition}

\begin{proof}
We first prove the continuous estimate at $q=Q$ for smooth
$G$-invariant functions. Integration over the group orbits gives
\[
\int_\Omega |u|^Q\,dx
=c_{m,n}\int_{\widehat\Omega}|u(s,t)|^Q
s^{m-1}t^{n-1}\,ds\,dt
\]
and
\[
\|u\|_{H^1(\Omega)}^2
=c_{m,n}\int_{\widehat\Omega}
\bigl(|u_s|^2+|u_t|^2+u^2\bigr)
s^{m-1}t^{n-1}\,ds\,dt,
\]
where $c_{m,n}:=|\mathbb S^{m-1}||\mathbb S^{n-1}|$.

Split the orbit domain into
\[
\widehat\Omega_s:=\widehat\Omega\cap\{s>t\},\qquad
\widehat\Omega_t:=\widehat\Omega\cap\{t>s\}.
\]
These two regions cover $\widehat\Omega$ up to a set of
measure zero. On the first region,
\[
\frac{R_-}{\sqrt2}<s<R_+,
\]
so the factor $s^{m-1}$ is bounded above and below by positive
constants. We may therefore remove this factor when estimating
the integrals on $\widehat\Omega_s$.

Write $t=|y|$, with $y\in\mathbb R^n$, and view the profile
as $u(s,|y|)$ on the bounded Lipschitz domain
\[
\left\{(s,y)\in\mathbb R\times\mathbb R^n:
s>|y|,\quad R_-^2<s^2+|y|^2<R_+^2\right\}.
\]
Polar integration in $y$ identifies the remaining weight
$t^{n-1}$ with the volume element in dimension $n$.
Thus the usual Sobolev inequality in dimension $n+1$ gives,
since $Q\leq2(n+1)/(n-1)$,
\[
\left(
\int_{\widehat\Omega_s}|u|^Q t^{n-1}\,ds\,dt
\right)^{2/Q}
\leq C\int_{\widehat\Omega_s}
\bigl(|u_s|^2+|u_t|^2+u^2\bigr)t^{n-1}\,ds\,dt.
\]
Restoring the bounded factor $s^{m-1}$ on both sides yields
\[
\begin{aligned}
\left(
\int_{\widehat\Omega_s}|u|^Q s^{m-1}t^{n-1}\,ds\,dt
\right)^{2/Q}
&\leq C\int_{\widehat\Omega_s}
\bigl(|u_s|^2+|u_t|^2+u^2\bigr)
s^{m-1}t^{n-1}\,ds\,dt\\
&\leq C\|u\|_{H^1(\Omega)}^2.
\end{aligned}
\]

On $\widehat\Omega_t$, we instead have
$R_-/\sqrt2<t<R_+$. Interchanging the two blocks removes
the bounded factor $t^{n-1}$ and leaves the weight $s^{m-1}$.
The same argument now uses the Sobolev inequality in dimension
$m+1$, whose critical exponent is exactly $Q$. Consequently,
\[
\left(
\int_{\widehat\Omega_t}|u|^Q s^{m-1}t^{n-1}\,ds\,dt
\right)^{2/Q}
\leq C\|u\|_{H^1(\Omega)}^2.
\]
Combining the estimates on the two regions and applying
Poincar\'e's inequality, we obtain
\[
\|u\|_Q\leq C\|u\|_{H^1(\Omega)}\leq C\|u\|.
\]
Smooth invariant functions are dense in $H$, as follows by
averaging $C_c^\infty(\Omega)$ approximations over $G$.
Hence this estimate holds for every $u\in H$.
Since $\Omega$ is bounded, H\"older's inequality also gives
\[
\|u\|_q\leq |\Omega|^{1/q-1/Q}\|u\|_Q
\leq C_q\|u\|,\qquad 1\leq q\leq Q.
\]

Finally, let $(u_j)$ be bounded in $H$. After passing to a
subsequence, Rellich compactness gives
\[
u_j\rightharpoonup u\quad\text{in }H,\qquad
u_j\to u\quad\text{in }L^2(\Omega).
\]
For $2<q<Q$, choose $\eta\in(0,1)$ such that
$1/q=\eta/2+(1-\eta)/Q$. Interpolation yields
\[
\|u_j-u\|_q
\leq\|u_j-u\|_2^\eta
\|u_j-u\|_Q^{1-\eta}\longrightarrow0,
\]
because the $L^Q$ norms are uniformly bounded. For $1\leq q<2$,
\[
\|u_j-u\|_q
\leq |\Omega|^{1/q-1/2}\|u_j-u\|_2\longrightarrow0.
\]
This proves compactness for every $1\leq q<Q$.
\end{proof}

Since $n\geq2$, we have $Q>2_N^*$. Thus the invariant space
allows a nonempty range of supercritical powers. We now impose
angular order on signed functions. Their two sign parts will
inherit opposite angular orders, which must be preserved in
the variational construction.

Define
\[
\C:=\{u\in H:u_\theta\geq0\text{ a.e.}\},
\]
\[
K^+:=\{v\in H:v\geq0,\ v_\theta\geq0\},\qquad
K^-:=\{v\in H:v\geq0,\ v_\theta\leq0\}.
\]
Thus $\C=\C_{m,n}$. We use the convention
$u^+:=\max\{u,0\}$ and $u^-:=\max\{-u,0\}$.

\begin{lemma}\label{lem:cone}
The set $\C$ is a closed convex cone in $H$. If $u\in\C$,
then $u^+\in K^+$ and $u^-\in K^-$. Conversely, if
$v\in K^+$, $w\in K^-$, and $vw=0$ a.e., then
$v-w\in\C$, with $(v-w)^+=v$ and $(v-w)^-=w$.
\end{lemma}

\begin{proof}
Linearity of the angular derivative gives the cone property
and convexity. Moreover,
\[
\|u_\theta\|_2\leq R_+\|u\|\qquad(u\in H).
\]
Hence if $u_j\in\C$ and $u_j\to u$ in $H$, then
$(u_j)_\theta\to u_\theta$ in $L^2$, so $u_\theta\geq0$.
Thus $\C$ is norm closed. Since it is convex, it is also
weakly closed.

The Sobolev chain rule gives
\[
(u^+)_\theta=\mathbf1_{\{u>0\}}u_\theta,
\qquad
(u^-)_\theta=-\mathbf1_{\{u<0\}}u_\theta,
\]
which proves the assertions about the sign parts.
Conversely, $(v-w)_\theta=v_\theta-w_\theta\geq0$,
so $v-w\in\C$. Nonnegativity and disjointness of $v,w$
identify them as its positive and negative parts.
\end{proof}

Throughout the next sections we assume
\begin{equation}\label{eq:p}
2<p<Q.
\end{equation}
The invariant embedding makes the energy
\[
J(u):=\frac12\|u\|^2
-\frac1p\int_\Omega a(r)|u|^p\,dx,
\qquad u\in H,
\]
a $C^1$ functional on $H$, with
\[
\langle J'(u),\varphi\rangle
=(u,\varphi)_H
-\int_\Omega a(r)|u|^{p-2}u\varphi\,dx
\qquad(u,\varphi\in H).
\]

We finish by specifying the weak formulation. In the
supercritical range, an arbitrary $H_0^1$ test function need
not belong to $L^p$. We therefore first use a definition in
which the nonlinear integral is well defined, and then show
that invariant solutions admit all $H_0^1$ test functions.

\begin{definition}\label{def:weak}
A weak solution of \eqref{eq:problem} is a function
$u\in H_0^1(\Omega)\cap L^p(\Omega)$ satisfying
\[
\int_\Omega\nabla u\cdot\nabla\varphi\,dx
=\int_\Omega a(r)|u|^{p-2}u\varphi\,dx
\qquad
\text{for every }\varphi\in H_0^1(\Omega)\cap L^p(\Omega).
\]
\end{definition}

We use the standard averaging argument for compact group
actions; see \cite[Section~5]{Palais1979}. In our supercritical
setting, we also verify the integrability needed for arbitrary
$H_0^1(\Omega)$ test functions. Let $\nu$ be the normalized Haar measure on $G$ and define
the group average by
\[
P_G\varphi(x):=\int_G\varphi(gx)\,d\nu(g).
\]
This operator maps $H_0^1(\Omega)$ into $H$ and is a
contraction in the Dirichlet norm and on each $L^q$ space.
For $\varphi\in H_0^1(\Omega)$, the invariant embedding gives
\[
\|P_G|\varphi|\|_p
\leq C\|P_G|\varphi|\|
\leq C\|\varphi\|_{H_0^1}.
\]
If $h\in L^{p'}(\Omega)$ is invariant, where
$p'=p/(p-1)$, Tonelli's theorem and orthogonal changes of
variables yield
\begin{equation}\label{eq:averaging-bound}
\int_\Omega|h\varphi|\,dx
=\int_\Omega|h|P_G|\varphi|\,dx
\leq C\|h\|_{p'}\|\varphi\|_{H_0^1}.
\end{equation}
In particular, the nonlinear pairing is integrable for every
$H_0^1$ test function when $h=a(r)|u|^{p-2}u$ and $u\in H$.

The same invariance gives
\[
\int_\Omega h\varphi\,dx
=\int_\Omega hP_G\varphi\,dx,
\]
and, for $u\in H$,
\[
\int_\Omega\nabla u\cdot\nabla\varphi\,dx
=\int_\Omega\nabla u\cdot\nabla(P_G\varphi)\,dx.
\]
Consequently, if $u\in H$ satisfies the weak equation for
tests in $H$, it satisfies it for every test in
$H_0^1(\Omega)$. This justifies the full weak formulations
used below.

\section{Attainment of the nodal Nehari level}
\label{sec:minimum}

We first construct a sign-changing function of least energy under
the angular-order constraint. The two signs must be treated
separately: the positive and negative parts will each satisfy a
Nehari identity, which fixes their amplitudes and prevents either
part from vanishing along a minimizing sequence. Compactness in
$H$ will then give attainment of the constrained minimum. Whether
this minimizer solves the full equation is the subject of the next
section.

Write
\[
B(u):=\int_\Omega a(r)|u|^p\,dx
\]
and define the nodal Nehari set by
\[
\M:=\left\{u\in\C:u^\pm\not\equiv0,\quad
\|u^+\|^2=B(u^+),\quad \|u^-\|^2=B(u^-)\right\}.
\]
Thus both sign parts lie on their respective Nehari constraints.
Our aim is to attain the level
\begin{equation}\label{eq:cnod}
c_{\C,\mathrm{nod}}:=\inf_{u\in\M}J(u),
\qquad \kappa_p:=\frac12-\frac1p=\frac{p-2}{2p}>0.
\end{equation}

The first observation is that every sign-changing function in
$\C$ can be placed on $\M$ by rescaling its two parts independently.
This is possible because the parts have disjoint supports, and
positive rescaling preserves their opposite angular orders.

\begin{lemma}\label{lem:scaling}
If $u\in\C$ and $u^\pm\not\equiv0$, the unique pair
$t_+,t_->0$ for which $t_+u^+-t_-u^-\in\M$ is
\[
t_\pm=\left(\frac{\|u^\pm\|^2}{B(u^\pm)}\right)^{1/(p-2)}.
\]
The function $(t,s)\mapsto J(tu^+-su^-)$ has its unique
global maximum on $(0,\infty)^2$ at $(t_+,t_-)$.
In particular, if $u\in\M$, this maximum occurs at $(1,1)$.
\end{lemma}

\begin{proof}
Lemma~\ref{lem:cone} shows that $tu^+-su^-\in\C$ for
every $t,s>0$. The sign parts are $tu^+$ and $su^-$, so their
Nehari identities are
\[
t^2\|u^+\|^2=t^pB(u^+),\qquad
s^2\|u^-\|^2=s^pB(u^-).
\]
Since $p>2$ and both parts are nonzero, these identities give
exactly the stated values of $t_+$ and $t_-$.

To identify the maximum, disjointness of the sign parts gives
\[
J(tu^+-su^-)=J(tu^+)+J(su^-).
\]
For every nonzero $v\in H$,
\[
\frac{d}{dt}J(tv)=t\|v\|^2-t^{p-1}B(v).
\]
This derivative is positive before
$t=(\|v\|^2/B(v))^{1/(p-2)}$ and negative afterwards.
Applying this observation to each part proves the claim.
\end{proof}

We now minimize on $\M$. The main point is to show that a weak
limit retains both signs. Once this is known, the preceding
rescaling restores any Nehari identities lost in the weak limit.
Minimality will force the rescaling factors to equal one and
will give strong convergence.

\begin{theorem}[Attainment of the nodal level]\label{thm:attainment}
The set $\M$ is nonempty and $0<c_{\C,\mathrm{nod}}<\infty$.
Every minimizing sequence has a subsequence converging strongly
in $H$ to a minimizer in $\M$.
\end{theorem}

\begin{proof}
Choose nonzero nonnegative smooth radial functions supported
in two disjoint concentric subannuli. Their difference belongs
to $\C$ and has two nonzero sign parts. Lemma~\ref{lem:scaling}
therefore produces an element of $\M$.

Next we obtain a lower bound for the size of each sign part.
By Proposition~\ref{prop:compactness} and the boundedness of $a$,
there is $C>0$ such that every nonzero $v\in H$ satisfying its
Nehari identity obeys
\[
\|v\|^2=B(v)\leq C\|v\|^p.
\]
Consequently,
\begin{equation}\label{eq:nonvanishing}
\|v\|^2\geq d:=C^{-2/(p-2)}>0.
\end{equation}
The same constant works for both sign parts. Moreover, on $\M$,
the two Nehari identities give
\begin{equation}\label{eq:Nehari-energy}
J(u)=\kappa_p\bigl(\|u^+\|^2+\|u^-\|^2\bigr)
\qquad(u\in\M).
\end{equation}
It follows that $c_{\C,\mathrm{nod}}\geq2\kappa_p d>0$.
Nonemptiness gives finiteness, and
\eqref{eq:Nehari-energy} shows that every minimizing sequence
$(u_j)\subset\M$ is bounded in $H$.

Let $(u_j)$ be such a sequence. After passing to a subsequence,
Proposition~\ref{prop:compactness} gives
\[
u_j\rightharpoonup u\quad\text{in }H,
\qquad u_j\to u\quad\text{in }L^p(\Omega)
\text{ and almost everywhere}.
\]
The cone is weakly closed, so $u\in\C$. The sign-part maps
are $1$-Lipschitz in $L^p$, and hence
$u_j^\pm\to u^\pm$ in $L^p$. They are also bounded in
$H_0^1(\Omega)$ and converge weakly there to $u^\pm$.
Thus
\[
B(u^\pm)=\lim_j B(u_j^\pm)
=\lim_j\|u_j^\pm\|^2\geq d.
\]
In particular, neither sign part disappears. Weak lower
semicontinuity further gives
\[
0<\|u^\pm\|^2\leq B(u^\pm).
\]

The limit may not yet satisfy equality in these two inequalities.
To recover the constraints, apply Lemma~\ref{lem:scaling} and set
\[
t_\pm:=\left(\frac{\|u^\pm\|^2}{B(u^\pm)}\right)^{1/(p-2)},
\qquad v:=t_+u^+-t_-u^-\in\M.
\]
We have $0<t_\pm\leq1$. Using minimality, the Nehari energy
identity for $v$, and the convergence above, we obtain
\begin{align*}
c_{\C,\mathrm{nod}}
&\leq J(v)
=\kappa_p\bigl(t_+^2\|u^+\|^2+t_-^2\|u^-\|^2\bigr)\\
&\leq\kappa_p\bigl(\|u^+\|^2+\|u^-\|^2\bigr)\\
&\leq\kappa_p\bigl(B(u^+)+B(u^-)\bigr)
=\lim_j J(u_j)=c_{\C,\mathrm{nod}}.
\end{align*}
Every inequality is therefore an equality. Since both sign parts
are nonzero, equality in the second line gives $t_+=t_-=1$.
Hence $u\in\M$ and $J(u)=c_{\C,\mathrm{nod}}$.

Finally, the Nehari identities and strong $L^p$ convergence give
\[
\|u_j\|^2
=B(u_j^+)+B(u_j^-)
\longrightarrow B(u^+)+B(u^-)=\|u\|^2.
\]
Weak convergence and convergence of the norms in the Hilbert
space $H$ imply $u_j\to u$ strongly in $H$.
\end{proof}

\section{The auxiliary operator and recovery of the equation}
\label{sec:criticality}

We have obtained a minimizer on the nodal Nehari set. Our next
task is to show that it solves the equation against arbitrary
test functions. The two Nehari identities alone only describe
variations of the positive and negative amplitudes; they do not
establish criticality under all variations. To recover the full equation, we introduce an auxiliary
Dirichlet operator and first prove that it preserves $\C$.

\subsection*{The cone-preserving auxiliary problem}

The angular-order property is first proved for a bounded source
with one weak derivative. These assumptions justify
differentiating the equation and determining the boundary values
of the angular derivative. We will then use truncation to treat
the nonlinear source, which need not be bounded initially.

\begin{lemma}\label{lem:linear-signed-cone}
Let $h\in H^1(\Omega)\cap L^\infty(\Omega)$ be invariant and
satisfy $h_\theta\geq0$ a.e. The weak solution of
$-\Delta v=h$ in $\Omega$, $v=0$ on $\partial\Omega$,
belongs to $\C$. No sign assumption on $h$ is required.
\end{lemma}

\begin{proof}
Lax--Milgram gives a unique solution in $H_0^1(\Omega)$.
Since the source and the domain are invariant, uniqueness gives
$v\in H$. Global Dirichlet estimates on the smooth annulus yield
$v\in W^{2,q}(\Omega)$ for every finite $q>1$, and hence
$v\in C^{1,\alpha}(\overline\Omega)$ for every $0<\alpha<1$;
see \cite{GilbargTrudinger}.

To prove angular monotonicity, set $z=v_\theta$. Away from the
orbit axes, differentiating the orbit equation gives
\begin{equation}\label{eq:angular-linear}
-\Delta_G z+Wz=h_\theta\geq0,\qquad
W:=\frac{m-1}{s^2}+\frac{n-1}{t^2},
\end{equation}
where
\[
\Delta_G:=\partial_{ss}+\partial_{tt}
+\frac{m-1}{s}\partial_s+\frac{n-1}{t}\partial_t.
\]
The differentiation is justified locally in distributions,
using $h\in H^1$ and $v\in W^{2,q}$. The coefficient $W$ is
nonnegative, which allows us to exclude negative values of $z$.

We first check its boundary values. The angular direction is
tangent to both spherical boundary components, where $v=0$;
thus $z=0$ there. Invariance and $C^1$ regularity give $v_s=0$
on $s=0$ and $v_t=0$ on $t=0$. Since $z=sv_t-tv_s$,
it also vanishes on the orbit axes.

For $\varepsilon>0$, let
$\zeta_\varepsilon:=(-z-\varepsilon)^+$.
Continuity and the zero boundary values show that its support
is compactly contained in the orbit domain and stays away from
the axes. In particular, the singular coefficients in
\eqref{eq:angular-linear} are bounded on its support, and
$-\zeta_\varepsilon$ is an admissible test function.
With $d\omega=s^{m-1}t^{n-1}\,ds\,dt$, testing gives
\[
\int_{\widehat\Omega}
\bigl(|\nabla\zeta_\varepsilon|^2
+W\zeta_\varepsilon^2
+\varepsilon W\zeta_\varepsilon\bigr)\,d\omega
=-\int_{\widehat\Omega}h_\theta\zeta_\varepsilon\,d\omega
\leq0.
\]
All terms on the left are nonnegative, so
$\zeta_\varepsilon=0$. Letting $\varepsilon\downarrow0$
proves $v_\theta\geq0$, as required.
\end{proof}

We now apply this property to the nonlinear source. Radiality
of $a$ and monotonicity of $t\mapsto |t|^{p-2}t$ preserve
angular order, including when $u$ changes sign. Truncation
allows us to use the preceding lemma without assuming
boundedness of $u$.

\begin{proposition}[Cone-preserving auxiliary operator]
\label{prop:auxiliary-operator}
For each $u\in H$, there is a unique $R(u)\in H$ such that
\begin{equation}\label{eq:R-weak}
\int_\Omega\nabla R(u)\cdot\nabla\varphi\,dx
=\int_\Omega a(r)|u|^{p-2}u\varphi\,dx
\qquad(\varphi\in H_0^1(\Omega)).
\end{equation}
The map $R:H\to H$ is continuous and satisfies
\[
\|R(u)\|\leq C\|u\|_p^{p-1},
\qquad R(\C)\subset\C.
\]
\end{proposition}

\begin{proof}
For $u\in H$, put $h_u:=a(r)|u|^{p-2}u\in L^{p'}(\Omega)$.
The invariant embedding gives
\[
\left|\int_\Omega h_u\varphi\,dx\right|
\leq\|h_u\|_{p'}\|\varphi\|_p
\leq C\|u\|_p^{p-1}\|\varphi\|
\qquad(\varphi\in H).
\]
Lax--Milgram on $H$ therefore gives a unique invariant solution
and the stated norm estimate. Averaging the test function as in
\eqref{eq:averaging-bound} extends this identity to every
$\varphi\in H_0^1(\Omega)$. Uniqueness in the full space follows
by testing the difference of two solutions with itself.

To obtain continuity, subtract the equations for $R(u)$ and
$R(v)$ and test with their difference. The invariant embedding
yields
\begin{equation}\label{eq:R-continuity}
\|R(u)-R(v)\|\leq C\|h_u-h_v\|_{p'}.
\end{equation}
The power map is continuous from $L^p$ to $L^{p'}$, so $R$
is continuous on $H$.

It remains to verify preservation of $\C$. For $u\in\C$, define
\[
T_k(t):=\max\{-k,\min\{t,k\}\},\qquad
h_k:=a(r)|T_k(u)|^{p-2}T_k(u).
\]
The Sobolev chain rule gives $h_k\in H^1(\Omega)\cap
L^\infty(\Omega)$. Since $a$ is radial,
\[
(h_k)_\theta
=a(r)(p-1)|T_k(u)|^{p-2}
\mathbf1_{\{|u|<k\}}u_\theta\geq0
\quad\text{a.e.}
\]
Lemma~\ref{lem:linear-signed-cone} implies that the Dirichlet
solution $v_k$ with source $h_k$ belongs to $\C$.
Moreover, $h_k\to h_u$ in $L^{p'}$ by dominated convergence.
The estimate used in \eqref{eq:R-continuity} gives
$v_k\to R(u)$ strongly in $H$. Closedness of $\C$ then
proves $R(u)\in\C$.
\end{proof}

\subsection*{An energy-decreasing deformation inside the cone}

Cone preservation provides an admissible direction: for
$u\in\C$, the whole segment joining $u$ to $R(u)$ remains
in $\C$. The following estimate shows that moving along this
segment decreases the energy strictly unless $u=R(u)$.
It is the energy estimate needed for the final argument.

\begin{lemma}\label{lem:resolvent-descent}
For $u\in\C$ and $0\leq\tau\leq1$, set
$u_\tau:=(1-\tau)u+\tau R(u)$. Then $u_\tau\in\C$ and
\begin{equation}\label{eq:descent}
J(u_\tau)\leq J(u)
-\tau\left(1-\frac\tau2\right)\|u-R(u)\|^2.
\end{equation}
In particular, the inequality is strict for $\tau>0$ whenever
$u\neq R(u)$.
\end{lemma}

\begin{proof}
Membership in $\C$ follows from convexity. Write
$\Phi(u):=B(u)/p$ and $d:=R(u)-u$. Since $a>0$ and
$p>2$, the functional $\Phi$ is convex. Therefore
\[
\Phi(u+\tau d)\geq\Phi(u)
+\tau\langle\Phi'(u),d\rangle.
\]
Combining this inequality with the quadratic expansion of the
Dirichlet term gives
\[
J(u+\tau d)\leq J(u)
+\tau\bigl[(u,d)_H-\langle\Phi'(u),d\rangle\bigr]
+\frac{\tau^2}{2}\|d\|^2.
\]
By \eqref{eq:R-weak},
$\langle\Phi'(u),d\rangle=(R(u),d)_H$. The expression in
brackets is consequently $-\|d\|^2$, proving
\eqref{eq:descent}.
\end{proof}

The segment above need not preserve the two Nehari identities.
To return to $\M$, we use the two amplitudes from
Lemma~\ref{lem:scaling}. The deformation will be turned off on
the boundary of a small parameter rectangle. This preserves the
opposite signs of the two Nehari defects on opposite sides,
allowing the Poincar\'e--Miranda theorem to recover both identities
at an interior point.

\subsection*{Recovery of the full weak equation}

\begin{theorem}\label{thm:nodal-solution}
Every minimizer of $J$ on $\M$ is a sign-changing
weak solution of \eqref{eq:problem}.
\end{theorem}

\begin{proof}
Let $u\in\M$ attain $c:=c_{\C,\mathrm{nod}}$. By
\eqref{eq:R-weak}, it suffices to prove $u=R(u)$. Suppose
otherwise. Fix $0<\varepsilon<1$ and consider
\[
D:=[1-\varepsilon,1+\varepsilon]^2,\qquad
g(t,s):=tu^+-su^-.
\]
Lemma~\ref{lem:scaling} shows that $g(D)\subset\C$ and that
\[
J(g(1,1))=J(u)=c,\qquad
J(g(t,s))<c\quad\text{for }(t,s)\neq(1,1).
\]
Thus $u$ is the unique maximum on this two-parameter family.
We will lower that maximum while keeping its boundary fixed.

Choose a continuous cutoff $\chi:D\to[0,1]$ with
$\chi=0$ on $\partial D$ and $\chi(1,1)=1$. Set
\[
L:=\max_{(t,s)\in D}\|R(g(t,s))-g(t,s)\|_p<\infty,
\qquad
m_0:=(1-\varepsilon)
\min\{\|u^+\|_p,\|u^-\|_p\}>0.
\]
Finiteness of $L$ follows from continuity of $R$ and compactness
of $D$. Choose $0<\tau_0<1$ so that $\tau_0L<m_0/2$ and define
\[
h(t,s):=g(t,s)+\tau_0\chi(t,s)
\bigl[R(g(t,s))-g(t,s)\bigr].
\]
By cone preservation and convexity, $h$ is a continuous map
from $D$ into $\C$, and $h=g$ on $\partial D$.

The smallness of $\tau_0$ ensures that neither sign disappears.
Indeed, $\|g(t,s)^\pm\|_p\geq m_0$ throughout $D$.
Using the $1$-Lipschitz property of sign parts in $L^p$, we obtain
\[
\|h(t,s)^\pm\|_p
\geq\|g(t,s)^\pm\|_p-\|h(t,s)-g(t,s)\|_p
\geq m_0-\tau_0L>\frac{m_0}{2}.
\]
Lemma~\ref{lem:resolvent-descent} also gives
$J(h(t,s))\leq J(g(t,s))$ everywhere. At $(1,1)$ the decrease
is strict because $\chi(1,1)=1$ and $u\neq R(u)$.
At every other point $J(g(t,s))<c$ already. Continuity on the
compact rectangle therefore implies
\begin{equation}\label{eq:lowered-maximum}
\max_{(t,s)\in D}J(h(t,s))<c.
\end{equation}

We now show that this lower-energy family still meets $\M$.
Define the two Nehari defects by the continuous map
\[
\mathcal F:H\to\R^2,\qquad
\mathcal F(z):=
\bigl(\|z^+\|^2-B(z^+),\ \|z^-\|^2-B(z^-)\bigr).
\]
Continuity follows from continuity of the sign-part maps in
$H_0^1$ and the embedding $H\hookrightarrow L^p$.
Since $h=g$ on $\partial D$ and $u\in\M$, on that boundary
we have
\[
\mathcal F(h(t,s))=\mathcal F(g(t,s))
=\bigl(\|u^+\|^2(t^2-t^p),\
\|u^-\|^2(s^2-s^p)\bigr).
\]
Because $p>2$, the first component is positive on
$t=1-\varepsilon$ and negative on $t=1+\varepsilon$.
The second component is positive on $s=1-\varepsilon$ and
negative on $s=1+\varepsilon$. The Poincar\'e--Miranda theorem
\cite{Kulpa}, applied to $\mathcal F\circ h$, yields a point
$(t_*,s_*)\in D$ such that
\[
\mathcal F(h(t_*,s_*))=0.
\]
Both sign parts of this function are nonzero, so
$h(t_*,s_*)\in\M$. But \eqref{eq:lowered-maximum} gives
$J(h(t_*,s_*))<c$, contradicting the definition of $c$.

Therefore $u=R(u)$. Equation~\eqref{eq:R-weak} is the full
weak formulation against every $H_0^1(\Omega)$ test function.
Since $u\in\M$, both signs are nonzero, completing the proof.
\end{proof}

We have now shown that every constrained minimizer is an actual
sign-changing solution on the whole annulus. The next section
establishes its regularity and uses minimality once more to prove
that it has exactly one positive and one negative nodal domain.

\section{Regularity and the number of nodal domains}
\label{sec:nodal-domains}

The preceding section shows that every minimizer on $\M$ is
a sign-changing weak solution. We now establish its classical
regularity and prove that it has exactly one positive and one
negative nodal domain. We begin with boundedness, which allows
us to apply standard elliptic regularity estimates.

\begin{proposition}\label{prop:nodal-regularity}
Every invariant weak solution $u\in H$ of \eqref{eq:problem}
belongs to $L^\infty(\Omega)$ and to
$C^{2,\alpha}(\overline\Omega)$ for some $\alpha\in(0,1)$.
\end{proposition}

\begin{proof}
We adapt the standard Moser iteration to the improved
invariant Sobolev embedding; see also the regularity
argument in \cite{CowanMoameni}. Choose $p<q<Q$, so that
$u\in L^q(\Omega)$. We apply the iteration separately to
$v=u^+$ and $v=u^-$.

For $L>0$ and $\ell\geq1$, set
\[
T_L(t):=\min\{t,L\}\quad(t\geq0),
\qquad w_L:=vT_L(v)^{\ell-1}.
\]
Test the equation with $vT_L(v)^{2\ell-2}$ for the positive
part and its negative for the negative part. These tests belong
to $H_0^1(\Omega)\cap L^p(\Omega)$. The chain rule and the
invariant Sobolev embedding give
\begin{equation}\label{eq:nodal-Moser-start}
\|w_L\|_q^2
\leq C\int_\Omega|\nabla w_L|^2\,dx
\leq C\ell^2\int_\Omega v^{p-2}w_L^2\,dx.
\end{equation}
Set
\[
d:=\frac{2q}{q-p+2},\qquad
\vartheta:=\frac{p-2}{q-2}\in(0,1).
\]
Since $2<d<q$, H\"older's inequality and interpolation yield
\[
\int_\Omega v^{p-2}w_L^2\,dx
\leq\|v\|_q^{p-2}\|w_L\|_d^2
\leq\|v\|_q^{p-2}
\|w_L\|_q^{2\vartheta}\|w_L\|_2^{2(1-\vartheta)}.
\]
Young's inequality absorbs the $L^q$ term in
\eqref{eq:nodal-Moser-start} and gives
\[
\|w_L\|_q^2\leq C\ell^b\|w_L\|_2^2,
\qquad b:=\frac{2}{1-\vartheta},
\]
where $C$ may depend on $\|v\|_q$, but not on $L$ or $\ell$.
Whenever $v\in L^{2\ell}(\Omega)$, letting $L\to\infty$
gives the integrability gain
\[
\|v\|_{q\ell}
\leq(C\ell^b)^{1/(2\ell)}\|v\|_{2\ell}.
\]
Starting with $\ell=1$ and iterating with
$\ell=(q/2)^j$, we obtain $v\in L^\infty(\Omega)$, since
\[
\sum_{j\geq0}\frac{1+j}{(q/2)^j}<\infty.
\]
Thus both sign parts, and hence $u$, are bounded.

The source $a(r)|u|^{p-2}u$ is now bounded. Global Dirichlet
estimates give $u\in W^{2,r}(\Omega)$ for every finite $r>1$,
and hence $u\in C^{1,\alpha}(\overline\Omega)$ for some
$\alpha\in(0,1)$. Since $a\in C^1$ and $p>2$, the source
is H\"older continuous. Schauder estimates then yield
$u\in C^{2,\alpha}(\overline\Omega)$, after decreasing
$\alpha$ if necessary; see \cite{GilbargTrudinger}.
\end{proof}

A nodal domain is a connected component of either $\{u>0\}$
or $\{u<0\}$. To count these components, we will retain one
of each sign and compare the resulting function with the
minimizer. The following lemma ensures that restricting to a
single component preserves the required angular order, without
assuming a smooth nodal boundary.

\begin{lemma}[Restriction to a nodal domain]
\label{lem:nodal-component-cone}
Let $u\in\C\cap C(\overline\Omega)$, with $u=0$ on
$\partial\Omega$. If $D$ is a connected component of $\{u>0\}$,
then $v_D:=u\mathbf1_D\in K^+$ and
\[
\nabla v_D=\mathbf1_D\nabla u\quad\text{a.e.}
\]
If $E$ is a connected component of $\{u<0\}$, then
$w_E:=(-u)\mathbf1_E\in K^-$ and
\[
\nabla w_E=-\mathbf1_E\nabla u\quad\text{a.e.}
\]
\end{lemma}

\begin{proof}
We treat a positive component $D$. We first justify that its
zero extension belongs to the Sobolev space. For $\varepsilon>0$,
the set
\[
\{x\in D:u(x)\geq\varepsilon\}
\]
is compactly contained in $D$, because $u$ vanishes on
$\partial D$. Choose a smooth cutoff supported in $D$ and
equal to one on this set. Multiplying $(u-\varepsilon)^+$
by the cutoff represents
\[
v_{D,\varepsilon}:=(u-\varepsilon)^+\mathbf1_D
\]
as an element of $H_0^1(\Omega)$, with
\[
\nabla v_{D,\varepsilon}
=\mathbf1_{D\cap\{u>\varepsilon\}}\nabla u.
\]
Dominated convergence for the functions and their gradients
as $\varepsilon\downarrow0$ gives $v_D\in H_0^1(\Omega)$
and the stated gradient identity.

We next verify symmetry and angular order. Since $m,n\geq2$,
every $G$-orbit is connected. The orbit of a point in $D$
lies in $\{u>0\}$ and therefore in the same component $D$.
Thus $D$ and $v_D$ are invariant.

The continuous orbit profile of $u$ is nondecreasing in
$\theta$: weak monotonicity gives this first for almost every
radius, and continuity extends it to every radius. If
$(r,\theta_1)$ represents a point of $D$, then for
$\theta\geq\theta_1$,
\[
u(r,\theta)\geq u(r,\theta_1)>0.
\]
The corresponding angular path remains in $\{u>0\}$ and hence
in $D$. Therefore, wherever $v_D$ is positive, it remains
nondecreasing toward larger angles. Together with $v_D\geq0$,
this proves that its whole profile is nondecreasing, so
$v_D\in K^+$.

For a negative component, apply the same restriction argument
to $-u$, using its nonincreasing angular profile and paths
toward smaller angles. This gives $w_E\in K^-$ and the
corresponding gradient identity.
\end{proof}

We can now use minimality. The full equation implies that each
retained component satisfies its own Nehari identity. Keeping
one component of each sign therefore gives an admissible
competitor, while any additional component contributes strictly
positive energy.

\begin{theorem}[Exactly two nodal domains]
\label{thm:two-nodal-domains}
Every minimizer of $J$ on $\M$ has exactly two nodal domains,
one positive and one negative.
\end{theorem}

\begin{proof}
Let $u$ be a minimizer. By Theorem~\ref{thm:nodal-solution}
and Proposition~\ref{prop:nodal-regularity}, it solves the full
equation and is continuous. Since $u\in\M$, both signs are
nonzero. Choose a positive component $D$ and a negative
component $E$, and set
\[
v:=u\mathbf1_D,\qquad w:=(-u)\mathbf1_E.
\]
Lemma~\ref{lem:nodal-component-cone} gives $v\in K^+$,
$w\in K^-$, and $vw=0$ a.e. In particular, $v-w\in\C$.
Testing the full weak equation with $v$ and $-w$ gives
\[
\|v\|^2=B(v),\qquad\|w\|^2=B(w).
\]
Hence $v-w\in\M$, and \eqref{eq:Nehari-energy} yields
\[
J(v-w)=\kappa_p\bigl(\|v\|^2+\|w\|^2\bigr).
\]

The gradient identities in the lemma also give
\[
\|u\|^2=\|v\|^2+\|w\|^2
+\int_{\Omega\setminus(D\cup E)}|\nabla u|^2\,dx.
\]
If another nodal domain existed, the restriction of $u$ to
that domain would be a nonzero element of $H_0^1(\Omega)$
and would have positive Dirichlet energy by Poincar\'e's
inequality. Therefore
\[
J(v-w)<\kappa_p\|u\|^2=J(u),
\]
contradicting minimality on $\M$. Thus $D$ and $E$ are the
only nodal domains.
\end{proof}

\begin{remark}
The nodal-domain count does not imply that the interior zero
set is a smooth connected hypersurface or a radial graph.
At an interior zero where $\nabla u\neq0$, the implicit
function theorem gives a locally regular hypersurface.
Singular nodal points are not excluded here.
\end{remark}
\section{Nonradiality and multiplicity}
\label{sec:nonradiality}

We now show when the minimizers constructed above are nonradial
and how different symmetries yield distinct solution pairs.
An estimate for thin annuli will make the required conditions
explicit.

\subsection*{Angular perturbations and nonradiality}

To exclude radial minimizers, we perturb the positive and negative
parts in opposite angular directions. Independent rescaling then
restores their Nehari identities, so the comparison remains
within $\M$. The Hardy condition below ensures that the resulting
competitor has lower energy.

\begin{theorem}\label{thm:nonradial-nodal}
If $(p-2)\beta_\Omega>2N$, every minimizer of $J$ on $\M$
is nonradial.
\end{theorem}

\begin{proof}
Suppose that a minimizer $u$ is radial. For
$\sigma\in\{+,-\}$, put
\[
A_\sigma:=\int_\Omega|\nabla u^\sigma|^2\,dx,
\qquad B_\sigma:=\int_\Omega\frac{(u^\sigma)^2}{|x|^2}\,dx.
\]
The Nehari identities give $B(u^\sigma)=A_\sigma$.

Choose the angular profile $Y$ to be the positive multiple of
$m\sin^2\theta-n\cos^2\theta$ normalized by
\[
\frac1{|\mathbb S^{N-1}|}
\int_{\mathbb S^{N-1}}Y^2\,dS=1.
\]
It is a degree-two spherical harmonic: before normalization,
it is the restriction of the harmonic polynomial
$m|x''|^2-n|x'|^2$ to the unit sphere. In particular,
\[
\int_{\mathbb S^{N-1}}Y\,dS=0,\qquad
-\Delta_{\mathbb S^{N-1}}Y=2NY,\qquad
Y_\theta\geq0.
\]
The zero mean removes the first-order energy terms, while
the angular monotonicity makes the perturbation admissible.

For small $|\varepsilon|$, define
\[
v_+(\varepsilon):=u^+(1+\varepsilon Y),\qquad
v_-(\varepsilon):=u^-(1-\varepsilon Y).
\]
The multipliers are positive, so both parts remain nonzero
and disjoint. Since $u^\pm$ are radial, for $\varepsilon\geq0$
we have
\[
(v_+)_\theta=\varepsilon u^+Y_\theta\geq0,
\qquad
(v_-)_\theta=-\varepsilon u^-Y_\theta\leq0.
\]
Thus $v_+(\varepsilon)\in K^+$ and
$v_-(\varepsilon)\in K^-$.

We next compute the effect of this perturbation on the two
terms in the Nehari identities. Set
\[
A_\sigma(\varepsilon):=\|v_\sigma(\varepsilon)\|^2,
\qquad
P_\sigma(\varepsilon):=B(v_\sigma(\varepsilon)).
\]
Radiality and the normalization of $Y$ give
\begin{align}
A_\sigma(\varepsilon)
&=A_\sigma+\varepsilon^2(A_\sigma+2NB_\sigma),
\label{eq:nodal-A-expansion}\\
P_\sigma(\varepsilon)
&=A_\sigma+\frac{p(p-1)}2A_\sigma\varepsilon^2
+o(\varepsilon^2).
\label{eq:nodal-P-expansion}
\end{align}
For the first identity, radial and tangential gradients are
orthogonal, and
$\int_{\mathbb S^{N-1}}|\nabla_{\mathbb S^{N-1}}Y|^2
=2N\int_{\mathbb S^{N-1}}Y^2$.
For the second, use radiality of $a$ and the uniform Taylor
expansion of $(1\pm\varepsilon Y)^p$.

To restore the constraints, let
\[
t_\sigma(\varepsilon):=
\left(\frac{A_\sigma(\varepsilon)}{P_\sigma(\varepsilon)}
\right)^{1/(p-2)},\qquad
U_\varepsilon:=t_+(\varepsilon)v_+(\varepsilon)
-t_-(\varepsilon)v_-(\varepsilon).
\]
Lemma~\ref{lem:scaling} and preservation of angular order under
positive rescaling give $U_\varepsilon\in\M$ for small
$\varepsilon\geq0$. Its energy is
\[
E(\varepsilon):=J(U_\varepsilon)
=\kappa_p\sum_{\sigma\in\{+,-\}}
A_\sigma(\varepsilon)^{p/(p-2)}
P_\sigma(\varepsilon)^{-2/(p-2)}.
\]
Expanding the scalar powers using
\eqref{eq:nodal-A-expansion}--\eqref{eq:nodal-P-expansion}
gives $E(0)=J(u)$, $E'(0)=0$, and
\[
E''(0)=-(p-2)(A_++A_-)+2N(B_++B_-).
\]
By the definition of $\beta_\Omega$,
$A_\sigma\geq\beta_\Omega B_\sigma$. Hence
\[
E''(0)
\leq\bigl[2N-(p-2)\beta_\Omega\bigr](B_++B_-)<0.
\]
It follows that $J(U_\varepsilon)<J(u)$ for small positive
$\varepsilon$. This contradicts minimality on $\M$ and proves
nonradiality.
\end{proof}

\begin{proof}[Proof of Theorem~\ref{thm:main}]
Attainment follows from Theorem~\ref{thm:attainment}, and
Theorem~\ref{thm:nodal-solution} gives the full weak equation.
Proposition~\ref{prop:nodal-regularity} and
Theorem~\ref{thm:two-nodal-domains} give regularity and the
nodal-domain count. Finally, Theorem~\ref{thm:nonradial-nodal}
gives nonradiality under \eqref{eq:nodal-Hardy-condition}.
\end{proof}

\subsection*{Different symmetries and multiplicity}

To obtain distinct solutions, we apply the construction to
different coordinate decompositions. For
$2\leq j\leq\lfloor N/2\rfloor$, let
\[
G_j:=O(N-j)\times O(j)
\]
act on the first $N-j$ and last $j$ standard coordinates.
The following observation ensures that solutions obtained from
two such groups cannot coincide unless they are radial.

\begin{lemma}\label{lem:separation}
A continuous function on $\Omega$ invariant under both $G_j$
and $G_\ell$, with $j<\ell$, is radial.
\end{lemma}

\begin{proof}
Write its two profiles as $F(s,t)$ and $H(\sigma,\tau)$,
using the decompositions $(N-j)+j$ and $(N-\ell)+\ell$.
Continuity allows the profiles to be evaluated on the orbit
axes. For any $\sigma,\tau\geq0$ with
$R_-^2<\sigma^2+\tau^2<R_+^2$, choose
\[
x:=\sigma e_1+\tau e_{N-\ell+1}.
\]
Its orbit variables for $G_\ell$ are $(\sigma,\tau)$.
Since $j<\ell$, both displayed coordinate directions belong
to the first $N-j$ block. Its orbit variables for $G_j$ are
therefore $(\sqrt{\sigma^2+\tau^2},0)$, giving
\[
H(\sigma,\tau)=F(\sqrt{\sigma^2+\tau^2},0).
\]
Every $G_\ell$-orbit has such a representative, so the function
depends only on $|x|$.
\end{proof}

\begin{proof}[Proof of Theorem~\ref{thm:main-multiplicity}]
For each $j\in\mathcal J$, Theorem~\ref{thm:main} applies
with $(m,n)=(N-j,j)$ and provides a nonradial solution $u_j$
with exactly two nodal domains. If $u_j=\pm u_\ell$ for
$j\neq\ell$, the same continuous function is invariant under
both groups. Lemma~\ref{lem:separation} would make it radial,
a contradiction. Thus the solutions are distinct up to sign.
Oddness of the equation gives the stated pairs $\{\pm u_j\}$.
\end{proof}

\subsection*{An explicit condition on thin annuli}

Both conclusions require a sufficiently large Hardy constant.
The next estimate expresses this requirement in terms of the
radii and shows that it holds on sufficiently thin annuli.

\begin{proposition}\label{prop:thin-annulus}
With $d:=R_+-R_-$, the Hardy constant satisfies
\begin{equation}\label{eq:thin-Hardy}
\beta_\Omega\geq\frac{\pi^2R_-^2}{d^2}
\left(\frac{R_-}{R_+}\right)^{N-1}.
\end{equation}
Consequently, for fixed $R_->0$,
$\beta_{A_{R_-,R_-+d}}\to\infty$ as $d\downarrow0$.
\end{proposition}

\begin{proof}
Let $v\in C_c^\infty(\Omega)$. Along each spherical direction
$\omega$, the function $r\mapsto v(r\omega)$ vanishes at the
two radial endpoints. Comparing the weight $r^{N-1}$ with
$R_-^{N-1}$ and $R_+^{N-1}$ and applying the one-dimensional
Dirichlet Poincar\'e inequality gives
\[
\int_{R_-}^{R_+}|\partial_r v|^2r^{N-1}\,dr
\geq\left(\frac{R_-}{R_+}\right)^{N-1}
\frac{\pi^2}{d^2}
\int_{R_-}^{R_+}|v|^2r^{N-1}\,dr.
\]
Integrating in $\omega$ and using
$|\nabla v|^2\geq|\partial_r v|^2$, we obtain the corresponding
bound for the full Dirichlet energy. Finally,
\[
\int_\Omega |x|^{-2}v^2\,dx
\leq R_-^{-2}\int_\Omega v^2\,dx
\]
gives \eqref{eq:thin-Hardy}. Density extends the estimate
to $H_0^1(\Omega)$.
\end{proof}

\begin{remark}
For a fixed decomposition, the interval
\[
2+\frac{2N}{\beta_\Omega}<p<Q_{m,n}
\]
is nonempty exactly when $\beta_\Omega>N(m-1)/2$.
For any nonempty index set $\mathcal J$ in
Theorem~\ref{thm:main-multiplicity}, each upper threshold
exceeds $2_N^*$. One may therefore first fix
\[
2_N^*<p<\min_{j\in\mathcal J}\frac{2(N-j+1)}{N-j-1}
\]
and then choose the annulus sufficiently thin to satisfy the
Hardy condition. For fixed $N$, this construction uses at most
$\lfloor N/2\rfloor-1$ standard coordinate decompositions.
\end{remark}

The following example gives two distinct solution pairs with
a common supercritical exponent.

\begin{example}\label{ex:two-pairs}
Take $N=6$, $p=16/5$, and $\Omega=A_{1,11/10}$.
For the decompositions $6=4+2$ and $6=3+3$,
\[
2_6^*=3<\frac{16}{5}<Q_{4,2}=\frac{10}{3}<Q_{3,3}=4.
\]
Estimate \eqref{eq:thin-Hardy} gives
\[
\beta_\Omega\geq100\pi^2(10/11)^5>10.
\]
Thus $(p-2)\beta_\Omega>12=2N$. For every coefficient
satisfying \eqref{eq:a}, Theorem~\ref{thm:main-multiplicity}
yields two distinct pairs of nonradial sign-changing solutions,
each with exactly two nodal domains, in a supercritical range.
\end{example}









\bibliographystyle{plain}
\bibliography{lsa}
\end{document}